\documentclass[12pt]{article}
\newcommand{\nocopyright}{
No Copyright\thanks{
The authors hereby waive all copyright
and related or neighboring rights to this work,
and dedicate it to the public domain.
This applies worldwide.
}}
\title{Laplace-isospectral hyperbolic 2-orbifolds are representation-equivalent}
\author{
Peter G. Doyle\thanks{Dartmouth College.}
\and
Juan Pablo Rossetti\thanks{
FaMAF-CIEM, Univ.\ Nac.\ C\'ordoba.
}}
\date{Version 2.0 dated 9 April 2014
\\ \nocopyright
}
\usepackage{amsthm,amssymb}
\usepackage{hyperref}
\usepackage{verbatim}
\usepackage{amsmath}
\usepackage{graphicx}
\usepackage{hyperref}

\newcommand{\floatplacement}{bthp}
\newcommand{\putfigwithsize}[4]{
\begin{figure}[\floatplacement]
\centerline{\mbox{\includegraphics[#1]{figures/#3}}}
\caption{#4}
\label{fig:#2}
\end{figure}
}
\newcommand{\putfig}[3]{
\begin{figure}[\floatplacement]
\centerline{\mbox{\includegraphics[height=6.0cm]{figures/#2.eps}}}
\caption{#3}
\label{fig:#1}
\end{figure}
}

\newcommand{\putfigpairverticalwithsize}[5]{
\begin{figure}[\floatplacement]
\centerline{\mbox{\includegraphics[#1]{figures/#3.eps}}}
\centerline{\mbox{\includegraphics[#1]{figures/#4.eps}}}
\caption{#5}
\label{fig:#2}
\end{figure}
}

\newcommand{\figref}[1]{\ref{fig:#1}}

\newtheorem{theorem}{Theorem}
\newtheorem{prop}{Proposition}
\newtheorem{corollary}{Corollary}
\newtheorem{lemma}{Lemma}
\newtheorem{fact}{Unproven Assertion}

\newcommand{\Lap}{\Delta}

\newcommand{\Z}{\mathbf{Z}}
\newcommand{\R}{\mathbf{R}}
\newcommand{\C}{\mathbf{C}}

\newcommand{\half}{\frac{1}{2}}

\newcommand{\linearspan}[1]{\langle #1 \rangle}

\newcommand{\fig}[1]{}

\newcommand{\Vol}{\mathrm{Vol}}
\newcommand{\Isom}{\mathrm{Isom}}
\newcommand{\choices}[1]{
\left\{
\begin{array}{ll}
#1
\end{array}
\right.
}
\newcommand{\goesinto}{\backslash}
\newcommand{\union}{\cup}
\newcommand{\card}{\#}
\newcommand{\Mbar}{{\bar{M}}}
\newcommand{\tr}{\mathrm{tr}}
\newcommand{\subgroup}{<}
\newcommand{\mathproofend}{\quad \qed}
\newcommand{\proofend}{$\quad \qed$}

\newcommand{\unqed}{\boxtimes}
\newcommand{\unproofend}{$\quad \unqed$}
\newcommand{\intover}{\int\limits}
\newcommand{\conj}[2]{\mathrm{Cl}(#1,#2)}
\newcommand{\cent}[2]{\mathrm{Z}(#1,#2)}
\begin{document}
\maketitle

\begin{abstract}
Let $M$ be a compact hyperbolic $2$-orbifold
(not necessarily connected).
We show that
the spectrum of the Laplacian on functions on $M$
determines the following:
the volume;
the total length of the mirror boundary;
the number (properly counted) of conepoints of each order;
and the number (properly counted) of primitive closed geodesics
of each length and orientability class.
This means that Laplace-isospectral hyperbolic 2-orbifolds are
representation-equivalent, and hence strongly isospectral.
\end{abstract}

\section{Statement}

Let $M$ be a compact hyperbolic 2-orbifold (not necessarily connected).
Denote the eigenvalues of the Laplacian acting on functions on $M$ by
\[
0=\lambda_0 \leq \lambda_1 \leq \ldots
.
\]
(Note that we don't necessarily have $0 < \lambda_1$
because we aren't assuming that $M$ is connected.)
We call the sequence $(\lambda_0,\lambda_1,\ldots)$
the \emph{Laplace spectrum} of $M$.
If two spaces have the same Laplace spectrum we call them
\emph{Laplace-isospectral}.
(Note that we don't simply call them `isospectral',
because this term is used in different ways by different authors.)

Our goal here will be to prove:

\begin{theorem}
\label{th:main}
Let $M$ be a compact hyperbolic $2$-orbifold (not necessarily connected).
The Laplace spectrum of $M$  determines, and is determined by,
the following data:
\begin{enumerate}
\item
the volume;
\item
the total length of the mirror boundary;
\item
the number of conepoints of each order, counting a mirror corner
as half a conepoint of the corresponding order;
\item
the number of closed geodesics of each length and orientability class,
counting a geodesic running along the boundary as half
orientation-preserving and half orientation-reversing,
and counting the $k$-fold iterate of a primitive geodesic
as worth $\frac{1}{k}$ of a primitive geodesic of the same length
and orientability.
\end{enumerate}
\end{theorem}

Of course the
Laplace spectrum determines other data as well, for example the number
of connected components.
The data we list here determine those other data, since they
determine the spectrum.

Theorem \ref{th:main} can be recast less picturesquely
as follows.
Associated to a 2-orbifold $M$ is a linear representation
$\rho_\Mbar$ of $\Isom(H^2)$
on functions on the frame bundle $\Mbar$ of $M$.
Associated to this representation is its character $\chi_\Mbar$,
a function on the set of conjugacy classes of $\Isom(H^2)$.
The geometrical data listed in Theorem \ref{th:main} are just
a way of describing geometrically the information conveyed
by the character $\chi_\Mbar$.
The character determines
the Laplace spectrum of $M$, and indeed, the spectrum of any other
natural operator.
Theorem \ref{th:main} tells us
that for hyperbolic $2$-orbifolds,
we can go back the other way, and recover the character from the 
Laplace spectrum.
From the character,
we get by general principles the linear equivalence class
of the representation $\rho_\Mbar$,
and the spectrum of any natural operator on any natural bundle over $M$.
Thus Laplace-isospectral hyperbolic $2$-orbifolds are
representation-equivalent and strongly isospectral.
We will discuss these matters further in section \ref{sec:imp} below;
for now we concentrate on Theorem \ref{th:main}.

\section{Examples}

In this section, we give examples of
the trade-offs between boundary and interior features that
are allowed for in the statement of Theorem \ref{th:main}.
Boundary corners on one orbifold may appear as conepoints on the other,
while boundary geodesics may migrate to the interior.
These examples should clarify just what it is we will be trying to prove,
and why we can't expect to prove something stronger.

The examples here will be obtained by glueing bunches of congruent
hyperbolic triangles.
The glueing patterns arise from transplantable pairs,
as described by 
Buser et al.\ \cite{bcds:drum}.

Examples of trading boundary and interior geodesics are very common.
(Indeed, practically any pair of isospectral hyperbolic 2-orbifolds
with relecting boundary exhibit this phenomenon.)
A variation on the famous example
of Gordon, Webb, and Wolpert
\cite{gordonWebbWolpert:drum}
yields a pair of planar hyperbolic 2-orbifolds
of types $*224236$ and $*224623$,
shown in Figure \figref{geotrade}.
This is likely the simplest pair of isospectral hyperbolic polygons,
both in having the smallest number of vertices and having the smallest area.
Each member of the pair is glued together from $7$ copies of
a prototype tile,
according to the glueing pattern labelled $7_3$ in  \cite{bcds:drum}.
The prototype tile is
a so-called \emph{$346$ triangle},
meaning a hyperbolic triangle with angles $\pi/3$, $\pi/4$, $\pi/6$.
\putfigwithsize{width=12cm}{geotrade}{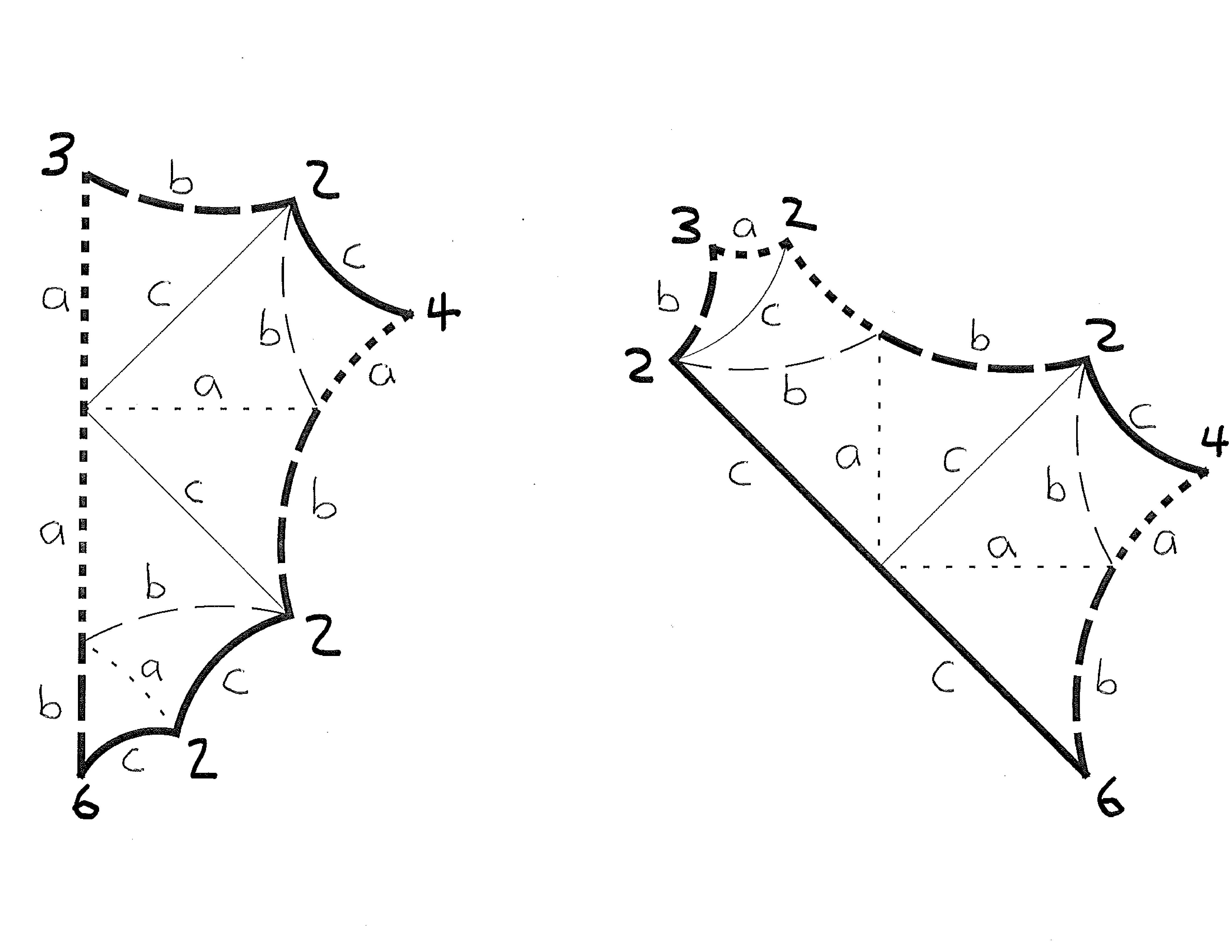}
{The pair $*224236$, $*224623$.}

There is no issue here with 
geodesics passing through the interior of any of the triangles that
make up these two orbifolds:
These can be matched so as to preserve
length, orientability, and index of imprimitivity.
But when it comes to geodesics that run along the edges of the triangles,
whether along the mirror boundary or in the interior of the orbifold,
it is necessary to balance boundary geodesics on one orbifold against
interior geodesics on the other, as provided for in Theorem \ref{th:main}.
To see this,
look at Figure \figref{geotrade},
and count geodesics on the two sides.
Beware that boundary geodesics turn back at corners of even
order,
but continue along around the boundary at corners of odd order.
Beware also of the way interior geodesics bounce when they hit the
boundary.
The resulting counts are indicated in Table \ref{table:count}.
The names `recto' and `verso' are short for `orientation-preserving'
and `orientation-reversing',
in analogy with the names for the front and back of a printed page.
The table only shows lengths for which there is at least one imprimitive
geodesic.
Trade-offs between boundary and interior geodesics continue at multiples
of these lengths.
\newcommand{\nl}{\\[0.2cm]}
\begin{table}[\floatplacement]
\centerline{Geodesics for $*224236$:}
\[
\begin{array}{c||cc|cc|cc}
\multicolumn{1}{c}{\ }&
\multicolumn{2}{c}{\mbox{boundary}}&
\multicolumn{2}{c}{\mbox{interior}}&
\multicolumn{2}{c}{\mbox{total}}
\nl
\mbox{length}&
\mbox{recto}&\mbox{verso}&
\mbox{recto}&\mbox{verso}&
\mbox{recto}&\mbox{verso}
\nl
\hline
2c&
\frac{3}{2}&\frac{3}{2}&
&&
\frac{3}{2}&\frac{3}{2}
\nl
2a+2b&
\half&\half&
1&1&
\frac{3}{2}&\frac{3}{2}
\nl
4c&
\frac{3}{2}\cdot\half&\frac{3}{2}\cdot\half&
1&&
\frac{7}{4}&\frac{3}{4}
\nl
4a+4b&
\half + \half \cdot \half& \half + \half \cdot \half&
2 \cdot \half&&
\frac{7}{4}&\frac{3}{4}
\nl
\end{array}
\]
\centerline{Geodesics for $*224623$:}
\[
\begin{array}{c||cc|cc|cc}
\multicolumn{1}{c}{\ }&
\multicolumn{2}{c}{\mbox{boundary}}&
\multicolumn{2}{c}{\mbox{interior}}&
\multicolumn{2}{c}{\mbox{total}}
\nl
\mbox{length}&
\mbox{recto}&\mbox{verso}&
\mbox{recto}&\mbox{verso}&
\mbox{recto}&\mbox{verso}
\nl
\hline
2c&
\half&\half&
1&1&
\frac{3}{2}&\frac{3}{2}
\nl
2a+2b&
\frac{3}{2}&\frac{3}{2}&
&&
\frac{3}{2}&\frac{3}{2}
\nl
4c&
\half + \half \cdot \half& \half + \half \cdot \half&
2 \cdot \half&&
\frac{7}{4}&\frac{3}{4}
\nl
4a+4b&
\frac{3}{2}\cdot\half&\frac{3}{2}\cdot\half&
1&&
\frac{7}{4}&\frac{3}{4}
\nl
\end{array}
\]
\caption{
\label{table:count}
Counting geodesics.
}
\end{table}

\newpage
Examples of trading corners for conepoints are not as common 
as examples of trading boundary and interior geodesics,
but there are still plenty.
Figures \figref{conetrade} and \figref{conetradehyper}
show how to construct a pair $6*2232233$, $23*22366$.
The transplantable pair used here is one of many found by John Conway
using the theory of quilts.
These will be described elsewhere, but there is no problem in verifying
transplantability here using the methods of  \cite{bcds:drum}.
Note how the order-6 conepoint moves to the boundary going one way, while the
order-2 and order-3 conepoints move to the boundary going the other way.
\putfig{conetrade}{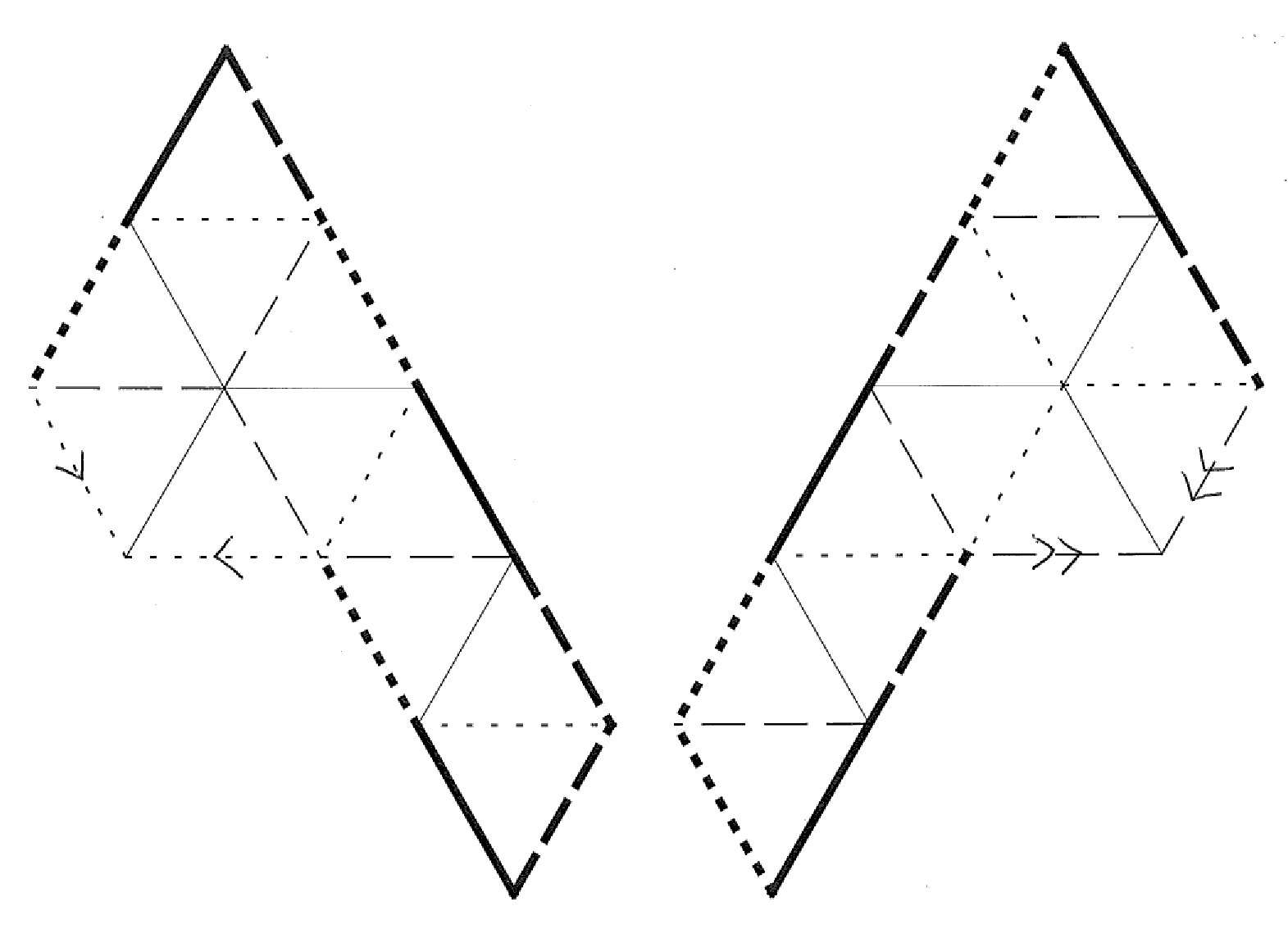}{
A transplantable pair of size 11.
}
\putfigwithsize{height=15cm}{conetradehyper}{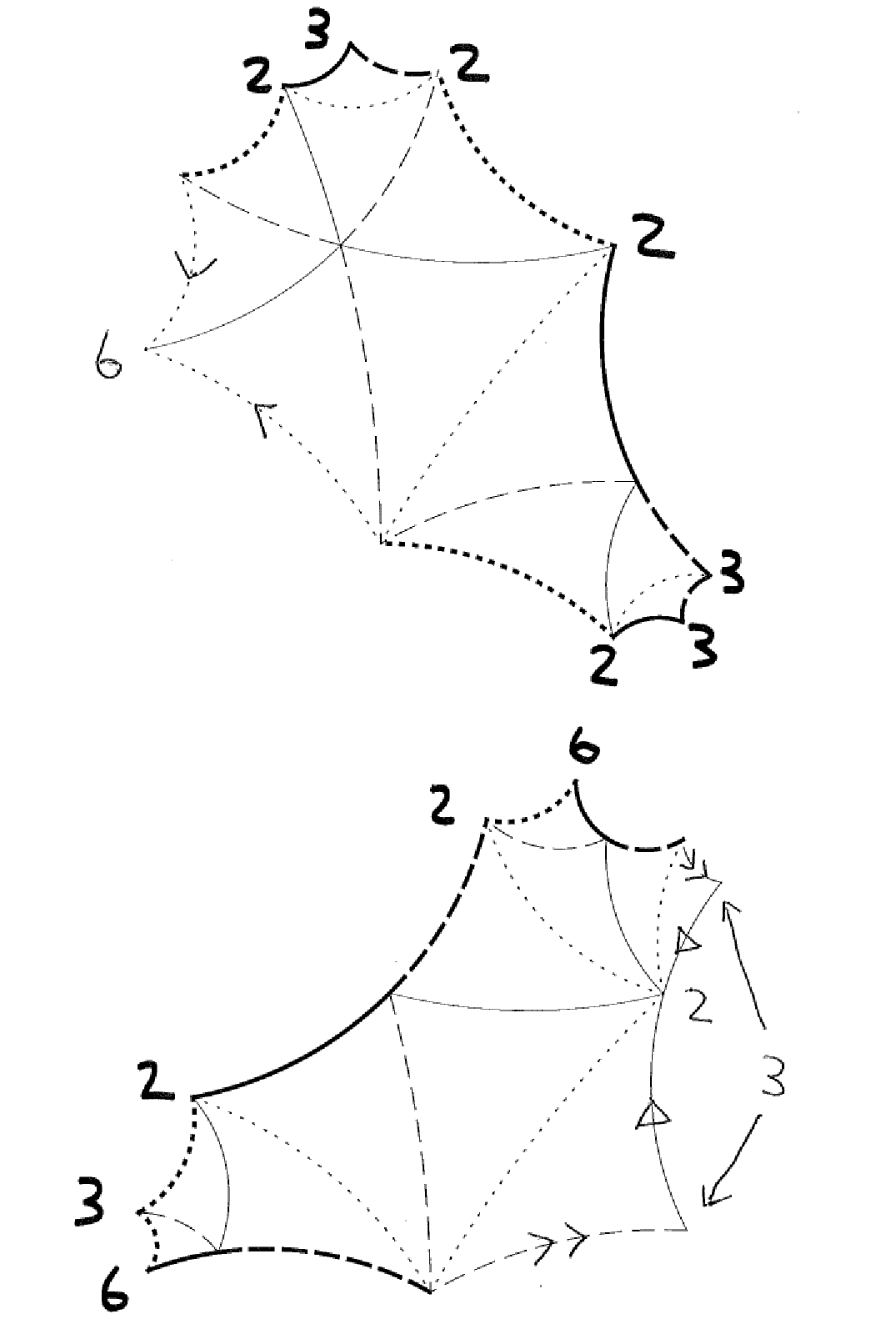}
{Replace the prototype equilateral triangle in Figure \figref{conetrade} by
a $366$ triangle
so that six 3-vertices come together in the left-hand diagram.
This yields a Laplace-isospectral pair of hyperbolic 2-orbifolds 
of types $6*2232233$ and $23*22366$.}

{\bf Note.}  The groups associated to the prototype triangles in these two examples
($346$ and $366$) both belong to the
finite list of \emph{arithmetic triangle groups}.
(Cf.\ Maclachlan and Reid
\cite{macreid}.)
The finite-index subgroup of a triangle group that belongs to an
orbifold obtained when the prototype triangle is arithmetic is
also arithmetic.
So our examples here are arithmetic.

\section{Background} 
Huber
\cite{huber:selbergI}
proved the result of Theorem \ref{th:main}
in the case of orientable hyperbolic 2-manifolds.
Huber used what would nowadays be seen as a version of the Selberg trace
formula,
which allows us to read off the lengths of geodesics from the spectrum
in a straight-forward way.
Doyle and Rossetti
\cite{doylerossetti:verso}
extended the result to non-orientable hyperbolic 2-manifolds.
In the non-orientable case
we can't simply read off the data about geodesics using the trace formula,
because of interference
between the spectral contributions of orientation-preserving
and orientation-reversing geodesics of the same length.
However, it turns out that
any possible scenario for matching spectral contributions would
require too many geodesics.

What we will show here is that,
as we indicated in
\cite{doylerossetti:verso},
it is a short step from non-orientable surfaces to general
orbifolds.
The reason is that the Selberg formula permits us to read off the data about
orbifold features, just as Huber read off the data about geodesics in
the case of orientable manifolds.
Then we have only to check that there is still no scenario
for matching the spectral contributions of the geodesics.

There are other ways to approach showing that the Laplace spectrum determines
the data about orbifold features.
By looking at the wave trace,
Dryden and Strohmaier
\cite{drydenstrohmaier}
proved the result of Theorem \ref{th:main}
for orientable hyperbolic 2-orbifolds.
While they did not consider non-orientable surfaces or orbifolds,
it seems likely that wave techniques could be used to
show that the Laplace spectrum determines all the orbifold data.
The argument would be essentially equivalent to the argument we give here,
only more complicated.

Another possible approach is via the heat equation.
By looking at short-time asymptotics of
the heat trace,
Dryden-Gordon-Greenwald-Webb  \cite{drydenetal} got information about the singular set of
general orbifolds
(for example, the volume of the reflecting boundary).
Their results 
yield information about orbifolds of variable curvature,
and in any dimension.
Restricted to hyperbolic $2$-orbifolds,
the results they state
don't yield complete information about the singular set.
All this information is there in the short-time asymptotics of the
heat trace, however,
and presumably it could be extracted using their approach,
by looking at higher and higher terms in the asymptotic expansion.

The beauty of the wave and heat approaches is that they can work in great
generality.
The Selberg method depends on having
spaces whose underlying geometry is homogeneous,
such as manifolds and orbifolds of constant curvature.
For studying such spaces,
it is a good bet that the Selberg method will
beat the wave and heat approaches.
Of course, the Selberg method can be used to treat the heat and wave
kernels;
the bet is that you will do better to consider a simpler kernel,
like the counting kernel that we use here.

Finally, just as it is true that hyperbolic manifolds whose geodesics match for
a sufficiently large number of lengths and twists
(in particular, all but a finite number)
must necessarily match for all lengths and twists
(cf. Kelmer \cite{kelmer:match}),
it is likely true on general principles that orbifolds with
matching geodesics must necessarily have matching orbifold features.
That would make Theorem \ref{th:main} follow almost immediately from the methods
used to prove the corresponding result 
\cite{doylerossetti:verso}
about $2$-manifolds, without explicitly computing the contributions
of orbifold features to the spectrum, as we do here.

\section{What we need from Selberg}

Here we assemble what we will need from Selberg for the proof
of Theorem \ref{th:main}.
All the ideas come from
Selberg
\cite{selberg:harmonic}.

Let $G=\Isom(H^2)$ be
the group of isometries
of $H^2$.
Note that $G$ has two components,
corresponding to orientation-preserving and orientation-reversing isometries.
A hyperbolic $2$-orbifold
can be written as a union of quotients of $H^2$ by discrete cocompact
subgroups $\Gamma_j \subset G$,
one for each connected component of $M$:
\[
M = \union_j \, \Gamma_j \goesinto H^2
.\]
We'll denote by $F(\Gamma_j)$ a fundamental domain for $\Gamma_j$.

In what follows, we could assume that $M$ is connected, and write
\[
M = \Gamma \goesinto H^2
,
\]
because the extension to the disconnected case presents no difficulties.
We prefer not to do this, in part to combat the common prejudice against
disconnected spaces.

Define the \emph{counting kernel} on $H^2$ by
\[
c(x,y;s)=\choices{1,&d(x,y)\leq s\\0,&d(x,y)>s}
.
\]
It tells when the hyperbolic distance $d(x,y)$ is at most $s$.
Define the \emph{counting trace}
\begin{eqnarray*}
C(s)
&=&
\sum_j \intover_{F(\Gamma_j)} \sum_{\gamma \in \Gamma_j} c(x,\gamma x;t) 
\,dx
\\&=&
\sum_j \intover_{F(\Gamma_j)} 
\card\{\gamma \in \Gamma_j:d(x,\gamma x)\leq s\}
\,dx
.
\end{eqnarray*}
The counting trace
tells (after dividing by the volume of $M$)
the average number of broken geodesic loops on $M$ of length at most $s$,

The reason for the name `counting trace' is that formally,
$C(s)$ is the trace of the linear operator
$L_s$ whose kernel is the counting kernel $c(x,y;s)$
pushed down to $M$.
$L_s$ associates to a function $f:M \to R$
the function whose value at $x \in M$
is the integral over a ball of radius $s$
of the lift of $f$ to the universal cover $H^2$.
$L_s$ is not a actually trace class operator, because of the discontinuity of
the counting kernel,
but our expression for $C(s)$ is still well-defined.

We will need the following two propositions,
which correspond to the two sides of the Selberg trace formula.

\begin{prop} \label{prop:spectrace}
The Laplace spectrum determines the counting trace.
\proofend
\end{prop}

Let $\conj{\gamma}{\Gamma}$ denote the conjugacy class
of $\gamma$ in $\Gamma$,
and let $\cent{\gamma}{\Gamma}$ denote the centralizer.
We have the usual one-to-one correspondence between
the set of cosets
$\cent{\gamma}{\Gamma} \goesinto \Gamma$
and the conjugacy class
$\conj{\gamma}{\Gamma}$,
where to the coset
$\cent{\gamma}{\Gamma} \delta$ we associate 
the conjugate $\delta^{-1} \gamma \delta$.

\begin{prop}
The counting trace can be expressed as a sum of contributions from
conjugacy classes of $\Gamma$:
\[
C(s) =
\sum_j 
\sum_{\conj{\gamma}{\Gamma_j}}
\Vol(\{x \in F(\cent{\gamma}{\Gamma_j}): d(x,\gamma x)\leq s \})
.
\mathproofend
\]
\end{prop}

The virtue of this proposition is that we can evaluate the summands
using simple hyperbolic trigonometry.

\section{Outline}

The counting trace $C(s)$
is built up of contributions from the conjugacy classes
of $\Isom(H^2)$.
The contribution from the identity is independent of $s$,
and equal to the volume of $M$.
The remaining conjugacy classes are of four kinds:
reflections, rotations, translations, and glide reflections.
Their contributions all vanish for $s=0$, which means that $C(0)$ is the
volume of $M$.
Translations and glide reflections make no contribution for $s$ smaller
than the length of the shortest closed geodesic,
so for small $s$ the only contributions are from the identity,
reflections, and rotations.

Begin by reading off the volume $C(0)$,
and subtracting this constant from $C(s)$.
In what is left,
reflections make a contribution of first order in $s$, proportional
to the total length of the mirror boundary.
Rotations contribute only to second order,
so we can read off the length of the mirror boundary, subtract out
its contribution, 
and what is left is (for small $s$)
entirely the result of rotations.

A simple computation now
shows that the contributions of
rotations through different angles are linearly independent.
This allows us to read off the orders of the conepoints,
and subtract out their contributions.
Once we've zapped the contributions of the identity,
reflections, and rotations,
what is left of the counting trace vanishes identically for small $s$,
and overall is entirely the result of
translations and glide reflections,
corresponding in the quotient to closed geodesics.

Here we run into the
main
difficulty in the proof,
which is that the contributions of orientation-preserving
and orientation-reversing geodesics of a given length
are not linearly independent.
Indeed, they are proportional, with a constant of proportionality
depending on the length of the geodesic.
So we can't simply read off these lengths.
Fortunately, we have already established that the spectrum determines
the lengths and twists of geodesics
for hyperbolic 2-manifolds,
(See Doyle and Rossetti \cite{doylerossetti:verso}.)
and the argument carries over directly to the orbifold case,
providing we take care to count geodesics along the boundary as
half orientation-preserving and half orientation-reversing.

\section{Details}

Denote by $R$ the combined length of all reflecting boundaries.
\begin{lemma}
The combined contribution of the reflecting boundaries to the
counting trace $C(s)$ is
\[
R \sinh \frac{s}{2}
.
\]
\end{lemma}

{\bf Proof.}
Let $\rho \in \Gamma_j$ be a reflection.
There are two possibilities for $\cent{\rho}{\Gamma_j} \goesinto H^2$:
a funnel or a planar domain.
(See Figure \figref{reflectorbs}.)
\putfig{reflectorbs}{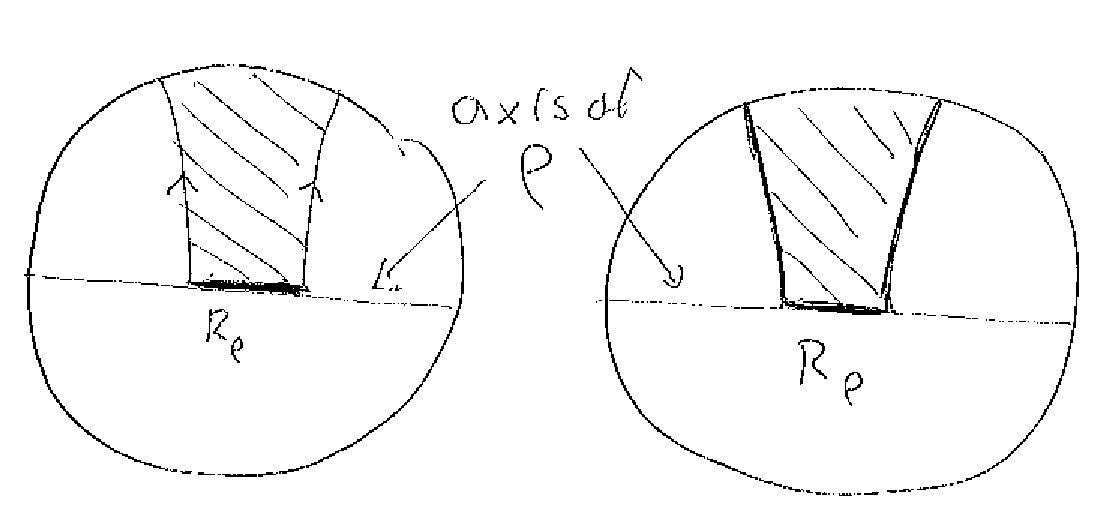}
{Possible quotients for the centralizer of a reflection.
The left shows a funnel, the right shows a planar domain.}
In either case, the contribution to $C(s)$ is $R_\rho \sinh \frac{s}{2}$,
where $R_\rho$ is the portion of $R$ attributable to $\conj{\rho}{\Gamma_j}$,
because this measures points in the quotient that are within $\frac{s}{2}$
of the axis of $\rho$, and hence within $s$ of their images under the
reflection $\rho$.
(See Figure \figref{reflectcontrib}.)
\putfig{reflectcontrib}{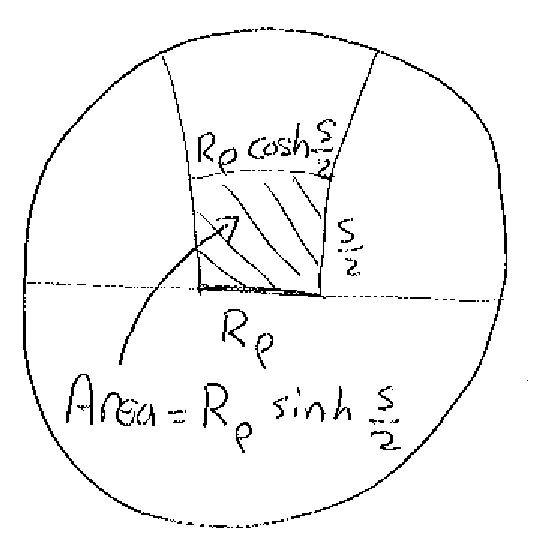}{The contribution of a reflection to the counting trace.}
Summing over the conjugacy classes of reflections in
all the groups $\Gamma_j$ gives
$R \sinh \frac{s}{2}$.
\proofend

Now we turn to conjugacy classes of rotations.
Each such is associated to a unique conepoint or boundary corner.
From a conepoint of order $n$ we get
$n-1$ rotations
through angles $\frac{2 \pi k}{n}$, $k=1,\ldots,n-1$.
For each such rotation $\gamma$ the centralizer $\cent{\gamma}{\Gamma_j}$
is cyclic of order $n$,
consisting of just these $n-1$ rotations, together with the identity.
The quotient $\cent{\gamma}{\Gamma_j} \goesinto H^2$ is
an infinite cone of cone angle $\frac{2 \pi}{n}$.
For a boundary corner of angle $\frac{\pi}{n}$
the centralizer is a dihedral group of order $2n$,
and the quotient is an infinite sector of angle $\frac{\pi}{n}$.
This quotient is half of a cone of angle $\frac{\pi}{n}$,
and it should be clear that the contribution to $C(s)$ of
rotations associated to a boundary corner is just half what it is for
a conepoint---which is why a boundary corner counts as half a cone-point.
So we can forget about boundary corners, and consider just conepoints.

Denote by $g_n(s)$ the total contribution to the 
counting trace of a conepoint of order $n$.
This contribution is the sum of the contributions of
the $n-1$ non-trivial conjugacy classes of rotations
associated to the conepoint.
Denote the contribution of the conjugacy class of the rotation through
angle $\frac{2 \pi k}{n}$
by $g_{k,n}(s)$.

\begin{lemma}
The total contribution of a conepoint of order $n$ to the counting trace
is $g_{n}(s) = \sum_{k=1}^n g_{k,n}(s)$,
where
\[
g_{k,n}(s) =
\frac{1}{n} f_{\frac{2 \pi k}{n}}(s),
\]
where
\[
f_\theta(s) =
2 \pi
\left (
\sqrt{1+\frac{\sinh^2 \frac{s}{2}}{\sin^2 \frac{\theta}{2}}} - 1
\right )
.
\]
\end{lemma}

{\bf Proof.}
It should be clear that
\[
g_{k,n}(s) =
\frac{1}{n} f_{\frac{2 \pi k}{n}}(s),
\]
where
$f_\theta(s)$ is the area $2 \pi(\cosh r - 1)$
of a hyperbolic disk of such a radius $r$ that
a chord of length $s$ subtends angle $\theta$.
(A fundamental domain for the cone hits a fraction $\frac{1}{n}$ of
any disk about the center of rotation.)
We just have to make sure that we have the correct formula for $f_\theta$.
The quantities $r, \theta, s$ satisfy
\[
\sinh r = \frac{\sinh \frac{s}{2}}{\sin\frac{\theta}{2}}
.
\]
(See Figure \figref{rtheta}.)
\putfig{rtheta}{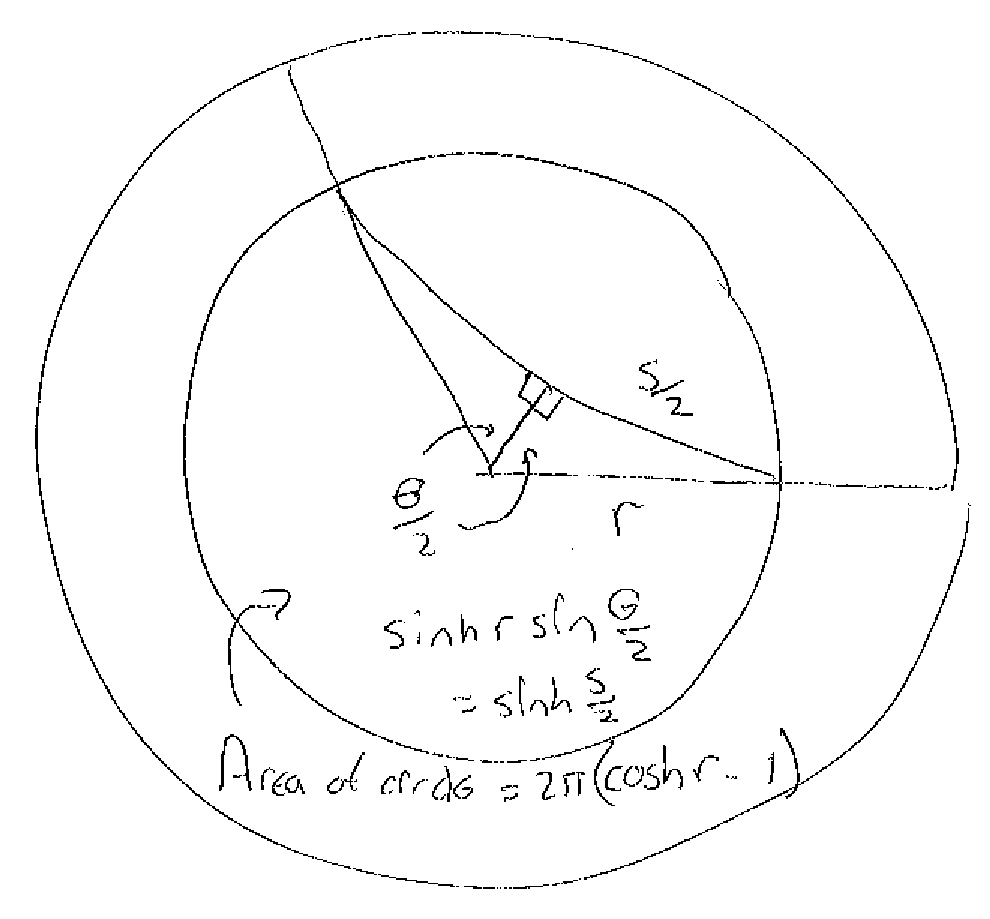}{The circle whose chord of length $s$ subtends angle $\theta$}
Thus
\begin{eqnarray*}
f_\theta(s)
&=& 
2 \pi  (\cosh r -1)
\\&=&
2 \pi 
\left (
\sqrt{1+\sinh^2 r}-1
\right )
\\&=&
2 \pi
\left (
\sqrt{1+\frac{\sinh^2 \frac{s}{2}}{\sin^2 \frac{\theta}{2}}} -1
\right )
.
\mathproofend
\end{eqnarray*}

Hopefully it will seem absolutely incredible that there might
be any non-trivial linear relation between these functions $f_\theta$.
There isn't:

\begin{lemma}
The functions $f_\theta, 0<\theta\leq \pi$ are linearly independent:
No non-trivial linear combination of any finite subset
of these functions vanishes identically on any interval $0 < s < S$.
\end{lemma}

{\bf Proof.}
It suffices to show linear independence of the family of functions
$\sqrt{1+a u}-1$, $1<a<\infty$, on all intervals $0 < u < U$.
(Set $u=\sinh^2 \frac{t}{2}$ and $a=\frac{1}{\sin^2 \frac{\theta}{2}}$.)
Actually, functions of this form remain independent when $a$ is allowed
to range throughout the punctured complex plane $\C - \{0\}$,
and not just over the specified
real interval.
A quick way to see this is to write in Taylor series
\begin{eqnarray*}
g(u)
&=&
\sqrt{1+u}-1
\\&=&
\frac{1}{1!}\left( \frac{1}{2} \right)
u
+
\frac{1}{2!}
\left( \frac{1}{2} \right)
\left( \frac{-1}{2} \right)
u^2
+
\frac{1}{3!}
\left( \frac{1}{2} \right)
\left( \frac{-1}{2} \right)
\left( \frac{-3}{2} \right)
u^3
+ \ldots
\\&=&
b_1 u + \frac{b_2}{2} u^2 + \frac{b_3}{6} u^3 + \ldots
,
\end{eqnarray*}
where all we will need is that
\[
b_1,b_2,\ldots \neq 0
.
\]
If a linear combination
$\sum_{i=1}^n c_i g(a_i u)$
is to vanish on $0 < u < U$, then all its derivatives
must vanish at $u=0$:
\[
\sum_{i=1}^n c_i a_i^k b_k = 0,
\quad k=1,2,\ldots
.
\]
(We ignore the constant term because all the functions involved here
vanish at $u=0$.)
Dividing by $b_k$ gives
\[
\sum_{i=1}^n c_i a_i^k = 0,
\quad k=1,2,\ldots
.
\]
Could this system of equations for $c_1,\ldots,c_n$
have a non-trivial solution?
If so, then the subsystem consisting of only the first $n$
of these equations would have a non-trivial solution:
\[
\sum_{i=1}^n c_i a_i^k = 0,
\quad k=1,\ldots,n
.
\]
Let $d_i=c_i a_i$, so that
\[
\sum_{i=1}^n d_i a_i^{k-1} = 0,
\quad k=1,\ldots,n
.
\]
The $n$-by-$n$ matrix $(a_i^{k-1})_{1 \leq i,j \leq n}$ of this system of linear equations is a
Vandermonde matrix, with determinant $\prod_{i<j} (a_i-a_j)$,
so as long as the $a_i$ are distinct the system has no non-trivial
solution.
\proofend

{\bf Note.} In appealing to the fact
that none of the derivatives of $\sqrt{1+u}$ happen to vanish
at $u=0$,
we are seizing upon an accidental feature of the problem.
A more robust proof can be based on the following lemma.

\begin{lemma}
Let $f(z)= c_k x^k$ be analytic at $0$.
Suppose $g(z) = \sum_{i=1}^n a_i f(\alpha_i z)$
with all $|\alpha_i|<1$.
If the Taylor series for $f$ and $g$ agree then $f$ is a polynomial.
\end{lemma}

{\bf Proof.}
We must have
\[
\sum_{i=1}^n a_i \alpha_i^k = 1
\]
for any $k$ for which $c_k \neq 0$.
If $f$ is not a polynomial,
we can choose $k$ big enough so that the left side is small,
and arrive at a contradiction.
\proofend

{\bf Note.}
It would be interesting to know the extent to which
this result still holds without the requirement that $a_i < 1$.
It will suffice to understand the case where all $|\alpha_i|=1$.
As an example, if we take
\[
f(z) = \frac{z+z^2}{1-z^3} = z + z^2 + z^4 + z^5 + z^7 + z^8 + \ldots
\]
we can write
\[
f(z) = -f(\omega z) - f(\omega^2 z)
\]
where
$\omega = e^{2 \pi i / 3}$
is a cube root of unity.
Are there examples where the $\alpha_i$ are not roots of unity?

{\bf Completion of the proof of Theorem \ref{th:main}.}
Subtract from $C(s)$ from the counting trace the volume $C(0)$
and the contributions of the reflectors.
Because of the linear independence of the functions $f_\theta$,
we can identify the number of conepoints of highest order.
(It should go without saying that this includes the half-conepoints
at boundary corners.)
Subtract out their combined contribution from $C(s)$.
Now we can identify the number of conepoints of the next highest order,
and subtract out their contribution.
Proceed until all conepoint contributions
have been removed.
What remains arises from the closed geodesics.

In \cite{doylerossetti:verso},
we showed that if $M$ is a manifold,
the characteristics of the closed geodesics are determined
by the counting trace,
and hence by the spectrum.
The reason is that
any possible scenario for a counterexample requires
the participation of too many geodesics.
Roughly speaking, any balance between primitive orientation-reversing and
orientation-preserving geodesics
gets destroyed when you look at imprimitive iterates of the geodesics,
so you have to keep adding new geodesics to compensate.
The proof
in \cite{doylerossetti:verso}
carries over here essentially without change.
\proofend

\section{Implications} \label{sec:imp}

Let $G=\Isom(H^2) \simeq PGL(2,\R)$ denote the group of isometries of $H^2$.
Note that $G$ has two connected components,
corresponding to orientation-preserving and orientation-reversing isometries.
A hyperbolic $2$-orbifold $M$
can be written as a union of quotients of $H^2$ by discrete cocompact
subgroups $\Gamma_j \subset G$,
one for each connected component of $M$:
\[
M = \union_j \, \Gamma_j \goesinto H^2
.\]
$M$ is naturally covered by
\[
\Mbar = \union_j \, \Gamma_j \goesinto G
.\]
$\Mbar$ is the bundle of all orthonormal frames of $M$, not just
those that have a particular orientation;
every point $x$ of $M$ is covered by two disjoint circles in $\Mbar$,
corresponding to the two orientation classes of frames of the tangent
space $T_x M$.
Note that there will be two connected components of $\Mbar$
for each orientable component of $M$, and one for each non-orientable
component,
because dragging a frame around an orientation-reversing loop
(e.g. one that simply bumps off a mirror boundary) will take you from
one orientation class of frames to the other.

$G$ acts naturally on the right on $\Mbar$,
and hence on $L^2(\Mbar)$.
This linear representation of $G$ is analogous to the matrix representation
of a finite permutation group.
A finite permutation representation $\rho$ is determined up to
linear equivalence by its character $\chi$,
with $\chi(g)=\tr \rho(g)$ counting the fixed points of $g$.
The situation here should be exactly analogous,
the only question being exactly how to define the character.
The answer comes from the Selberg trace formula.

For now let's extend the discussion to an arbitrary unimodular Lie group
$G$, possibly disconnected.
Let $\Gamma \subgroup G$ be a discrete subgroup with compact quotient
$\Mbar = \Gamma \goesinto G$.
(Notice that we make no mention here of $M$, though in the intended
application $\Mbar$ will arise from a homogeneous quotient
$M = \Gamma \goesinto G / K$.)
Denote by $\conj{g}{G}$ the conjugacy class of $g$ in $G$,
and by $\cent{g}{G}$ the centralizer of $g$ in $G$.
Introduce Haar measure $\rho^g$ on $Z_G(g)$,
normalized in a consistent (i.e., $G$-equivariant) way.
Attribute to $\conj{\gamma}{\Gamma}$ the \emph{weight}
$\rho^\gamma(\cent{\gamma}{\Gamma} \goesinto \cent{\gamma}{G})$,
and define the \emph{character} $\chi_\Mbar$ to be the function associating
to $g \in G$ the total weight of all conjugacy classes $\conj{\gamma}{\Gamma}$
for which $\conj{\gamma}{G} = \conj{g}{G}$.
Extend the definition of $\chi_\Mbar$
to $\Mbar = \union_j \, \Gamma_j \goesinto G$ by linearity.

As in the finite case, the character is a \emph{central function},
which means that $\chi_\Mbar(gh)=\chi_\Mbar(hg)$ for all $g,h \in G$,
or what is the same,
$\chi_\Mbar(g)$ depends only on the conjugacy class $\conj{g}{G}$.
In analogy to the finite case,
we can think of $\chi_\Mbar(g)$
as measuring (in appropriate units) the size of the
fixed point set
$\{x \in \Mbar: x g = x \}$.

\begin{prop}[Berard \cite{berard:repequiv}]
The character $\chi_\Mbar$ determines the representation $L^2(\Mbar)$
up to linear equivalence.
\end{prop}

{\bf Proof.}
Using a version of Selberg's trace formula
we can write the trace of the integral operator associated to any
smooth function of compact support on $G$
in terms of the values of
the character $\chi_\Mbar$.
(Cf.\ Wallach
\cite[Theorem 2.1]{wallach:stf},
Selberg \cite[(2.10) on p.\ 66]{selberg:harmonic}.)
By standard representation theory,
these traces determine the representation.
For details, see B\'erard
\cite{berard:repequiv}.
\proofend

\begin{prop}[DeTurck-Gordon \cite{deturckgordon:selberg}]
The character $\chi_\Mbar$ determines the trace of any natural
operator on any natural vector bundle over $M$.
\end{prop}
{\bf Proof.}
Selberg again.
\proofend

With this preparation, we have the following corollaries of
Theorem \ref{th:main}.

\begin{corollary}[Character is determined]
The Laplace spectrum of a hyperbolic 2-orbifold $M$ (not necessarily connected)
determines
the character $\chi_\Mbar$.
\proofend
\end{corollary}

\begin{corollary}[Representation-equivalence]
Laplace-isospectral hyperbolic $2$-orbifolds
(not necessarily connected)
determine
equivalent representations of $\mathrm{Isom}(H^2)$.
\proofend
\end{corollary}

\begin{corollary}[Strong isospectrality]
If two compact hyperbolic 2-orbifolds
(not necessarily connected)
are Laplace-isospectral 
then they have the same spectrum for any natural operator acting on
sections of any natural bundle.
\proofend
\end{corollary}

\section{Counterexamples} \label{sec:counterex}

The analog of Theorem \ref{th:main} fails for flat 2-orbifolds.
In the connected case it does go through more or less by accident,
just because there aren't many connected flat 2-orbifolds.
But there are counterexamples among disconnected flat 2-orbifolds.
We described such examples briefly in
\cite{doylerossetti:verso}.

\newcommand{\premier}{{$\frac{1}{2}+\frac{1}{6}=\frac{2}{3}$}}
Here's our premier example, which we call \emph{the {\premier} example}.
Let $H_1$ denote the standard hexagonal flat torus
$\linearspan{(1,0),(-\frac{1}{2},\frac{\sqrt{3}}{2})} \goesinto \R^2$.
$H_1$ has as 2-, 3-, and 6-fold quotients
a $2222$ orbifold $H_2$ (this is a regular tetrahedron); a $333$
orbifold $H_3$; and a $236$ orbifold $H_6$.
(See Figures \figref{hexfig} and \figref{hexorb}.)
\putfigwithsize{height=14cm}{hexfig}{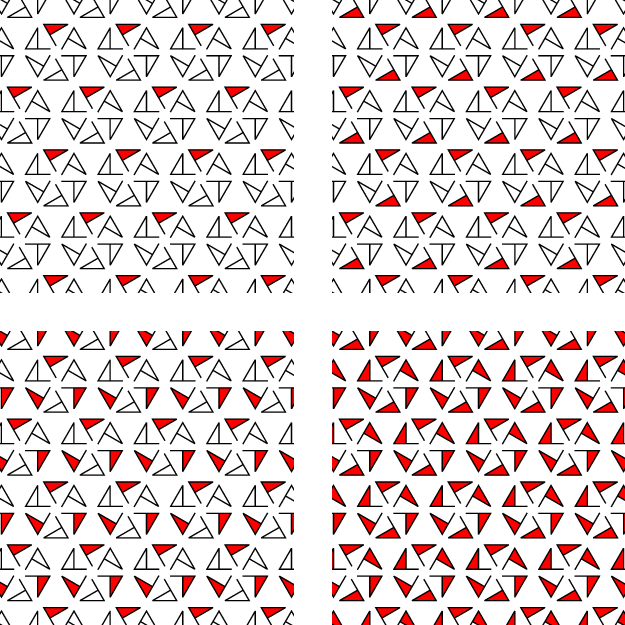}{The universal covers of the standard hexagonal
torus $H_1$ and its quotient orbifolds $H_2$, $H_3$, and $H_6$.}
\putfigwithsize{height=14cm}{hexorb}{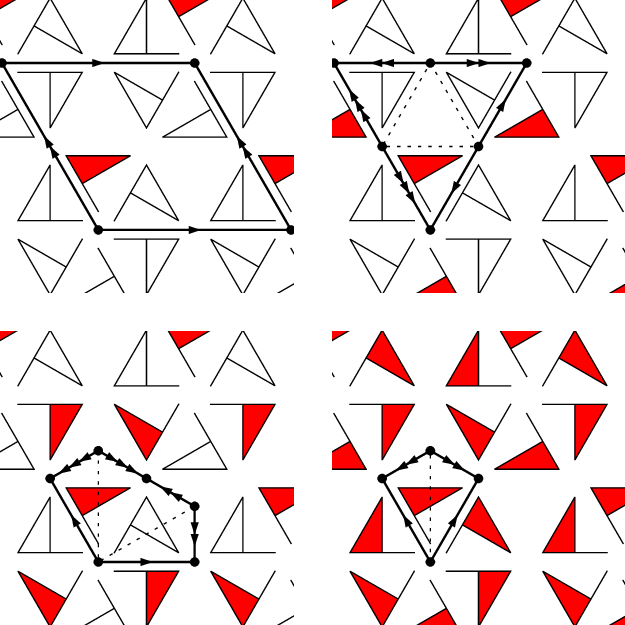}{The orbifolds
$H_1 = \bigcirc$, $H_2 =2222$, $H_3=333$, and $H_6=236$.}
Spectrally,
\begin{equation}
H_2 + H_6 \equiv 2 H_3
,
\end{equation}
meaning that these spaces have the same Laplace spectrum.
(In fact, as we'll discuss below, they are isospectral
for the Laplacian on $k$-forms for $k=0,1,2$.)
Of course these two spaces match as to volume
($\frac{1}{2}+\frac{1}{6}=\frac{2}{3}$)
and number of components
($1+1=2$).
But conepoints do not match.
On the left we have a $2222$ and a $236$, so
five conepoints of order $2$;
one of order $3$;
one of order $6$.
On the right we have two $333$s, so nine conepoints of order $3$.

The reason this example is possible is that in the flat case,
the contributions of rotations to the counting kernel differ
only by a multiplicative constant.
Lumping together the contribution of all the rotations associated with a
conepoint of order $n$, we get something proportional to
$\frac{n^2-1}{n}$.
(This nice simple formula was discovered by Dryden-Gordon-Greenwald \cite{drydenetal}.)
For a general 2-orbifold with variable curvature, it measures
the contribution of conepoints to the short-time asymptotics of the heat trace.
In the flat case, the short-time asymptotics determine the
entire contribution, because the contributions
of all rotations are exactly proportional.

Combining contributions of all conepoints,
we get the following totals:
\[
\begin{array}{rcl}
\mbox{orbifold}&\mbox{contribution}&\mbox{examples}
\\
\hline
\\
2222:& 4 \cdot \frac{3}{2}=6
&H_2,T_2
\\
333:& 3 \cdot \frac{8}{3}=8
&H_3
\\
244:& \frac{3}{2} + 2 \cdot \frac{15}{4} = 9
&T_4
\\
236:& \frac{3}{2}+\frac{8}{3}+\frac{35}{6}=10
&H_6
\end{array}
\]
The isospectrality $H_2+H_6 \equiv 2H_3$ arises because $6+10=2\cdot8$ makes conepoint
contributions match; we've already observed that the volumes match; and geodesic contributions match because on both sides the covering manifolds are $2H_1$.

We get a second isospectrality
\[
H_1 + H_3 + H_6  \equiv 3 H_2
\]
because $8+10=3\cdot6$, $1+\frac{1}{3}+\frac{1}{6}=3 \cdot \frac{1}{2}$,
and on both sides the covering manifolds are $3H_1$.

Combining these two relations yields other Laplace-isospectral pairs, e.g.:
\[
H_1 + 3 H_3 \equiv 4 H_2;
\]
\[
H_1 + 4 H_6 \equiv 5 H_3;
\]
\[
2 H_1 + 3 H_6 \equiv 5 H_2
.
\]
By linearity, all these relations satisfy the condition of
having equal volume, number
of components, and contributions from conepoints.
But look here:
We're talking, in effect, about formal combinations of
$H_1,H_2,H_3,H_6$, so we're in a space of dimension $4$.
We have three linear conditions, but to our surprise,
the subspace they determine has dimension $2$.
Our three linear conditions are not independent.
If we match volume and number of components,
agreement of conepoint contributions follows for free.
But why?  We don't know.

A similar coincidence happens for square tori.
Let $T_1$ denote the standard square torus $\Z^2 \goesinto \R^2$.
$T_1$ has as 2-, and 4-fold quotients
a $2222$ orbifold $T_2$ 
and a $244$ orbifold $T_4$.
We're in a 3-dimensional space,
so we expect to be out of luck when we impose 3 constraints,
but in fact we have the relation
\[
T_1+2T_4 \equiv 3T_2
.\]
We check equality of volume $1+\frac{2}{4} = \frac{3}{2}$
and number of components $1+2=3$, and then we find that equality of conepoint
contributions follows for free: $2\cdot 9 = 3 \cdot 6$.
Again, why?

These relations apply only to the spectrum of the Laplacian on functions,
not $1$-forms.
When we apply the Selberg formula to $1$-forms,
there is still interference between spectral contributions of conepoints,
which allows our {\premier} example to continue to slip through:
It is isospectral for $1$-forms as well as for functions.
But we get additional constraints, which taken with those coming from
the $0$-spectrum force there to be the same number of torus components
on each side,
and this wipes out the other examples.
The reason we have to have the same number of torus components is that this
is just half the dimension of the space of `constant' or `harmonic'
$1$-forms, meaning those belonging to the $0$-eigenspace of the Laplacian.
And since we are expecting just one more linear constraint, this must
be it.
So, no additional mystery here,
though when you do the computation it still seems a little surprising.

The isospectrality of the {\premier} example
does not hold for all natural operators.
because the left side admits harmonic differentials
that the right side does not.
Specifically,
look at the real part of a (non-holomorphic) quadratic differential,
i.e. an expression which in any conformally
correct chart with coordinate $z$ takes the form
$\Re f(z) dz^2$,
where $f$ is an arbitrary (not necessarily holomorphic) complex-valued
function.
We can rewrite this as
\[
g(x,y) (dx^2-dy^2) + h(x,y) (2 dx dy)
\]
where $g$ and $h$ are real.
Such quadratic differentials are sections of a natural flat bundle,
and in a flat coordinate system the Laplacian passes through to act on
$g$ and $h$.
Wherever the section doesn't vanish we can set
\[
Re f(z) dz^2 = 0
\]
and solve for $dz$ to get an unoriented line element.
A harmonic section (i.e. an eigensection belonging to eigenvalue 0)
will have $g$ and $h$ constant in any flat coordinate system,
and the line elements will be parallel.
In the {\premier} example only
the orbifold $H_2$ admits a globally parallel line
element like this, so the spectra for the Laplacian on these sections
are different.

\section{Higher dimensions}
Version 1 of this paper
\cite{doylerossetti:orbv1}
contains, among many other things omitted here, the following conjecture:
`It seems likely that in any dimension, if two hyperbolic
(or spherical) orbifolds are isospectral on $k$-forms for all $k$,
then they are representation-equivalent.'

Lauret, Miatello, and Rossetti
\cite{lmr:lens}
have shown that this is false in the spherical case.
They give many examples of lens spaces that are all-forms-isospectral
but not representation-equivalent, beginning in dimension 5.
(For a simpler approach to these examples, see DeFord and Doyle
\cite{hodge}.)

In the hyperbolic case, the question remains open.

As for the flat case,
we saw in Section \ref{sec:counterex} an example of disconnected $2$-orbifolds
that are all-forms-isospectral but not representation-equivalent.
We have been hoping to
find higher-dimensional examples that are connected, and perhaps even
manifolds and not just orbifolds,
but so far we've had no success with this.

\bibliography{orb}
\bibliographystyle{hplain}
\end{document}